\documentclass[11pt]{amsart}
\usepackage{amsmath, amssymb,amsfonts,amscd}

\title[Diffeomorphism groups of the real line]{The diffeomorphism groups of the real line are pairwise bihomeomorphic}

\author[T. Banakh]{Taras Banakh}
\address[T. Banakh]{Department of Mathematics,  Ivan Franko National University of Lviv, Ukraine, and \newline
Instytut Matematyki, Uniwersytet Humanistyczno-Przyrodniczy Jana Kocanowskiego w Kielcach, Poland}  \email{t.o.tbanakh@gmail.com}

\author[T. Yagasaki]{Tatsuhiko Yagasaki}
\address[T. Yagasaki]{Division of Mathematics,  Graduate School of Science and Technology, Kyoto Institute of Technology, Kyoto, 606-8585, Japan} \email{yagasaki@kit.ac.jp} 

\subjclass[2000]{57S05, 58D05, 58D15}
\keywords{The diffeomorphism group, the Whitney topology, the box product.}

\newcommand{\IR}{\mathbb R}

\newcommand{\IZ}{\mathbb Z}
\newcommand{\ID}{\mathsf D}

\newcommand{\II}{\mathbb I}

\newcommand{\w}{\omega}

\newcommand{\E}{\mathsf E}

\newcommand{\D}{\mathcal D}

\newcommand{\supp}{\operatorname{supp}}
\newcommand{\cbox}{\boxdot}

\newcommand{\boldprod}{\mbox{\large$\bold\Pi$}}
\newcommand{\boldsum}{\mbox{\large$\bold\Sigma$}}
\newcommand{\ls}{\leqslant}

\newcommand{\I}{\mathcal{I}}

\newcommand{\HH}{\mathcal H}

\newcommand{\id}{\mathrm{id}}

\newcommand{\cl}{\mathrm{cl}}

\newcommand{\bs}{\boldsymbol} 
\newcommand{\ds}{\displaystyle}

\newtheorem{theorem}{Theorem}[section]
\newtheorem{proposition}{Proposition}[section]
\newtheorem{lemma}{Lemma}[section]
\newtheorem{corollary}{Corollary}[section]
\newtheorem{claim}{Claim}[section]

\theoremstyle{definition}

\newtheorem{problem}{Problem}[section]

\begin{document}

\baselineskip 5mm 

\begin{abstract} For an $r = 0, 1, \cdots, \infty$, by $\D^r(\IR)$, $\D^r_+(\IR)$, $\D^r_c(\IR)$ we denote respectively the groups of $C^r$ diffeomorphisms, orientation-preserving $C^r$ diffeomorphisms, and compactly supported $C^r$ diffeomorphisms of the real line.
We think of these groups as bitopologies spaces endowed with the compact-open $C^r$ topology and the Whitney $C^r$ topology. We prove that all the triples $(\D^r(\IR),\D^r_+(\IR),\D^r_c(\IR))$, $0 \le r\le\infty$, are pairwise bitopologically equivalent, which allows to apply known results on the topological structure of homeomorphism groups of the real line to recognizing the topological structure of the diffeomorphisms groups of $\IR$.
\end{abstract}

\maketitle 

\section{Introduction} 

Groups of homeomorphisms and diffeomorphisms of smooth manifolds are classical objects in Differential Geometry and Topology, see \cite{Hirsch}. In this paper we shall study the compact-open and Whitney $C^r$ topologies on the group $\D^r(\IR)$ of $C^r$ diffeomorphisms of the real line for 
$r=0, 1, \cdots, \infty$. The group $\D^r(\IR)$ will be considered together with the subgroup $\D_+^r(\IR)\subset\D^r(\IR)$ of orientation-preserving $C^r$ diffeomorphisms of $\IR$, and the subgroup $\D_c^r(\IR)\subset\D_+(\IR)$  of $C^r$ diffeomorphisms $h:\IR\to\IR$ that have compact support
$\supp(h)=\cl_\IR\{x\in\IR:h(x)\ne x\}$.
In such a way, we obtain the triple of diffeomorphism groups $(\D^r(\IR),\D^r_+(\IR),\D^r_c(\IR))$. 


The group $\D^r(\IR)$ carries two natural topologies: the compact-open $C^r$ topology and the Whitney $C^r$ topology, see \cite{Hirsch}, \cite{Illman}. It will be convenient to consider those two topologies on $\D^r(\IR)$ simultaneously by thinking of $\D^r(\IR)$ as a bitopological space.

By a bitopological space we understand a triple $(X,w,s)$ consisting of a set $X$ and two topologies $w\subset s$ on $X$, called the {\em weak} and {\em strong topologies} of the bitopological space.  A function $f:X\to Y$ between bitopological spaces is called {\em weakly} (resp. {\em strongly}) {\em continuous} if $f$ is continuous with respect to the weak (resp. strong) topologies on $X$ and $Y$. A function $f:X\to Y$ is {\em bicontinuous} if $f$ is both weakly and strongly continuous.
Two bitopological spaces are called {\em bihomeomorphic} (resp. {\em weakly homeomorphic}, {\em strongly homeomorphic}) if there is a bijective map $f:X\to Y$ such that both $f$ and $f^{-1}$ are bicontinuous (resp. weakly continuous, strongly continuous).

Examples of bitopological spaces very often arise in Topology and Functional Analysis. 
For example, for two topological spaces $X,Y$ the space $C(X,Y)$ of continuous functions  from $X$ to $Y$ can be seen as a bitopological space whose weak and strong topologies are the compact-open and Whitney topologies. We recall that the Whitney topology on $C(X,Y)$ is generated by the base consisting of the sets $U_\Gamma=\{f\in C(X,Y):\{(x,f(x)):x\in X\}\subset U\}$ where $U$ runs over open subsets of $X\times Y$.

Another example of a bitopological space is the Cartesian product $\boldprod_{i\in\I}X_i$ of topological spaces $X_i$, $i\in\I$. Its weak topology is the Tychonov product topology while the strong topology is the box-product topology. The latter is  generated by the base consisting of products $\prod_{i\in\I}U_i$ of open subsets $U_i\subset X_i$, $i\in\I$. The product $\prod_{i\in\I}X_i$ endowed with the box-product topology  is denoted by $\square_{i\in\I}X_i$. If each space $X_i$, $i\in\I$, has a distinguished point $*_i\in X_i$, then the product $\boldprod_{i\in\I}X_i$ contains an important subspace
$$\boldsum_{i\in\I}X_i=\{(x_i)_{i\in\I}\in\prod_{i\in\I}X_i:|\{i\in\I:x_i\ne*_i\}<\aleph_0\}$$called the {\em $\sigma$-product} of the pointed spaces $X_i$. 


If all the spaces $X_i$, $i\in\I$, coincide with some fixed space $X$, then the spaces  $\boldprod_{i\in\I}X$ and $\boldsum_{i\in\I}X$ are denoted by $\boldprod^\I X$ and $\boldsum^\I X$, respectively. It is interesting to remark that the bitopological space $\boldprod^\I X$ is bihomeomorphic to the function space $C(\I,X)$ where $\I$ is endowed with the discrete topology.

The spaces $\boldprod^\I X$ and $\boldsum^\I X$ endowed with the Tychonov product topology are denoted by $\mbox{\large $\Pi$}^\I X$ and $\mbox{\large$\Sigma$}^\I X$.
The same spaces endowed with the box-product topology are denoted by $\square^\I X$ and $\cbox^\I X$, respectively. So, in a sense, the bitopological space $\boldprod^\I X$ decomposes into two topological spaces: $\mbox{\large$\Pi$}^\I X$ and $\square^\I X$.
The same remark concerns the bitopological space $\boldsum^\I X$ that decomposes into two topological spaces $\mbox{\large$\Sigma$}^\I X$ and $\cbox^\I X$. 

Bitopological spaces that are (locally) homeomorphic to $\boldsum^\w\IR$ were characterized in the paper \cite{BS} that also presents many examples of such spaces appearing in Topology, Topological Angebra or Analysis. 
In particular, it is proved in \cite{BS} that a bitopological linear space $X$ is bihomeomorphic to $\boldsum^\w\IR$ if and only if $X$ has infinite algebraic dimension, the weak topology of $X$ is metrizable and $X$ endowed with its strong topology is the direct limit of finite-dimensional compacta, see \cite[2.2]{BS}. 
\smallskip

Another natural example of a bitopological space whose bitopological structure has been recognized is the function space $C(\IR)$ endowed with the compact-open and Whitney topologies. According to the (proof of the) main result in \cite{GZ1} the function space $C(\IR)$ is bihomeomorphic to the countable power $\boldprod^\w l_2$ of the separable Hilbert space $l_2$. The same result is true for a bitopological space $C^r(\IR)$ of $C^r$-differentiable functions on the real line. The weak and strong topologies of $C^r(\IR)$ are induced from the weak and strong topologies of the product $\boldprod_{m\in\IZ}C^r(\II_m)$ via the embedding $$C^r(\IR)\hookrightarrow\boldprod_{m\in\IZ}C^r(\II_m),\;\;f\mapsto (f|\II_m)_{m\in\IZ}.$$
Here $C^r(\II_m)$ is the Fr\'echet space of all $C^r$ differentiable functions $f:\II_m\to\IR$ on the closed interval $\II_m=[m,m{+}1]$. The (metrizable) topology of $C^r(\II_m)$ is generated by the seminorms $$\|f\|_{n}=\max_{0 \le k\le n}\max_{x\in\II_m}|f^{(k)}(x)| \ \ \ (n \in \w, n \le r).$$

The weak topology of the bitopological space $C^r(\IR)$ is referred to as the  {\em compact-open $C^r$ topology} while the strong topology as the {\em Whitney $C^r$ topology}, see \cite{Hirsch}, \cite{Illman}. 

In the space $C^r(\IR)$ consider the subspace $$C^r_{c}(\IR)=\{f\in C^r(\IR):f|\IR\setminus K\equiv 0\mbox{ for some compact subset $K\subset\IR$}\}$$consisting of $C^r$ differentiable functions with compact support.

The following theorem describing the bitopological structure of the pair $(C^r(\IR),C_c^r(\IR))$ was implicitly proved in \cite{GZ1}, \cite{GZ2}.

\begin{theorem}[Guran, Zarichnyi] For every $r \in \w \cup \{ \infty \}$ the pair $(C^r(\IR),C^r_{c}(\IR))$ is bihomeomorphic to the pair $(\boldprod^\w l_2,\boldsum^\w l_2)$. 
\end{theorem}

We say that for bitopological spaces $X,Y$ and their subsets $X'\subset X$ and $Y'\subset Y$ the pairs $(X,X')$ and $(Y,Y')$ are {\em bihomeomorphic} if there is a bihomeomorphism $h:X\to Y$ such that $h(X')=Y'$. 
\smallskip

In this paper we concentrate at studying the bitopological structure of the diffeomorphism group $\D^r(\IR)\subset C^r(\IR)$ and its subgroups $\D^r_+(\IR)$ and $\D^r_c(\IR)$.
For $r=0$ the group $\D^0(\IR)$ turns into the homeomorphism group $\HH(\IR)$ of the real line. In this case the triple $(\D^0(\IR),\D^0_+(\IR),\D^0_c(\IR))$ is denoted by $(\HH(\IR),\HH_+(\IR),\HH_c(\IR))$. In \cite{BMS} it was proved that the pair $(\HH_+(\IR),\HH_c(\IR))$ is weakly  and strongly homeomorphic to the pair $(\boldprod^\w l_2,\boldsum^\w l_2)$  (surprisingly, but we do not know if these two pairs are bihomeomorphic!)

In this paper we shall extend the mentioned result to the pairs $(\D^r_+(\IR),\D^r_c(\IR))$ for all $r\le\infty$. 
Our principal result is:

\begin{theorem}\label{main} For every $r \in \w \cup \{ \infty \}$ the triple $(\D^r(\IR),\D^r_+(\IR),\D^r_c(\IR))$ is bihomeomorphic to the triple $(\HH(\IR),\HH_+(\IR),\HH_c(\IR))$. 
\end{theorem}

We say that for bitopological spaces $X,Y$ and their subsets $X_1,X_2\subset X$ and $Y_1,Y_2\subset Y$ the {\em triples} $(X,X_1,X_2)$ and $(Y,Y_1,Y_2)$ are {\em bihomeomorphic} if there is a bihomeomorphism $h:X\to Y$ such that $h(X_1)=Y_1$ and $h(X_2)=Y_2$. 

By \cite{BMS} the pairs of homeomorphism groups $(\HH_+(\IR),\HH_c(\IR))$ and $(\HH(\IR),\HH_c(\IR))$ endowed with the Whitney topology are homeomorphic to the pair  $(\square^\w l_2,\cbox^\w l_2)$. This fact combined with Theorem~\ref{main} implies

\begin{corollary}\label{cor1} For every $r \in \w \cup \{ \infty \}$ the pairs $(\D^r(\IR),\D^r_c(\IR))$ and $(\D^r_+(\IR),\D^r_c(\IR))$ of diffeomorphism groups endowed with the Whitney $C^r$ topology are homeomorphic to the pair $(\square^\w l_2,\cbox^\w l_2)$.
\end{corollary}

On the other hand, the pair of homeomorphism groups $(\HH_+(\IR),\HH_c(\IR))$ endowed with the compact-open topology is homeomorphic to the pair $(\mbox{\large$\Pi$}^\w l_2,\mbox{\large$\Sigma$}^\w l_2)$, see \cite{BMS}. Combining this fact with Theorem~\ref{main} we get another 

\begin{corollary}\label{cor2} For every $r \in \w \cup \{ \infty \}$ the pair 
$(\D^r_+(\IR),\D^r_c(\IR))$ of diffeomorphism groups endowed with the compact-open $C^r$ topology is homeomorphic to the pair $(\mbox{\large$\Pi$}^\w l_2,\mbox{\large$\Sigma$}^\w l_2)$.
\end{corollary}

\section{Preliminaries}

In the sequel $n,m,k,i,s$ will denote elements of the set $\w$ of all non-negative integers. On the other hand, $r\in\w\cup\{\infty\}$ can also accept the infinite value $\infty$. Each number $n\in\w$ will be identified with the discrete topological space $\bold n=\{0,\dots,n-1\}$. By $\IR_+=(0,\infty)$ we denote the open half-line.

For an $r\in\w\cup\{\infty\}$ by $C^r(\IR)$ we denote the space of all $C^r$ differentiable functions $f:\IR\to\IR$. This space has a rich algebraic structure. First, it is an algebra with respect to the addition $f+g$ and the multiplication $fg$ of functions $f,g\in C^r(\IR)$. On the other hand, $C^r(\IR)$ is a monoid with respect to the operation $f\circ g$ of composition of functions. The identity map $\id_\IR$ of $\IR$ is the two-sided unit of this monoid, while $\D^r(\IR)$ is the subgroup of invertible element of $C^r(\IR)$. The group $\D^r_+(\IR)$ of orientation-preserving $C^r$ diffeomorphisms of $\IR$ lies in the submonoid $C^r_+(\IR)\subset C^r(\IR)$ consisting of functions $f:\IR\to\IR$ with strictly positive derivative $f'$. If $r=0$ then $C^r_+(\IR)$ is the space of strictly increasing continuous functions on $\IR$. 
The classical theorem of Inverse Function implies that the group $\D^r_+(\IR)$ coincides with the set of surjective functions $f\in C^r_+(\IR)$.

For a subset $B\subset\IR$ we put $C^r(\IR,B)=\{f\in C^r(\IR):f(\IR)\subset B\}$. If $r=0$, then we write $C(\IR)$ instead of $C^0(\IR)$.
Functions $f\in C^\infty(\IR)$ will be called {\em smooth}.

For each function $f\in C^r(\IR)$ and every $n \in \w$, $n \leq r$ we can calculate the (finite or infinite) norm
$$\|f\|_n=\max_{0 \le k\le n}\sup_{x\in\IR}|f^{(k)}(x)|.$$

If no other topology on $C^r(\IR)$ is specified, then we think of $C^r(\IR)$ as a bitopological space endowed with the compact open and Whitney $C^r$ topologies. Those are induced from the Tychonov product and box-product topologies of $\boldprod_{n\in\IZ}C^r(\II_n)$ under the embedding 
$$e:C^r(\IR)\to \boldprod_{n\in\IZ}C^r(\II_n), \ e:f\mapsto(f|\II_n)_{n\in\IZ}.$$
 Here $C^r(\II_n)$ is the Fr\'echet space of $C^r$ differentiable functions on the interval $\II_n=[n,n{+}1]$, endowed with the topology generated by the seminorms $\ds \|f\|_m=\max_{0 \le k\le m}\max_{x\in\II_n}|f^{(k)}(x)|$ for $m \in \w$, $m \le r$.
\smallskip

\section{The differential operators $\partial^D_A$}

Given a point $a\in\IR$ and a number $n\in\w$, consider the differential operator
$$\partial^{(n)}_a:C^\infty(\IR)\to \IR,\quad\partial^{(n)}_a:f\mapsto f^{(n)}(a)$$assigning to a $C^\infty$ differentiable function its $n$-th derivative at the point $a$.

For a subset $D\subset\w$, the operators $\partial^{(n)}_a$, $n\in D$, compose an operator
$$\partial^D_a:C^\infty(\IR)\to \IR^D,\quad\partial^D_a:f\mapsto(f^{(n)}(a))^{n\in D}.$$
The power $\IR^D$ is considered as a bitopological space whose weak and strong topologies coincide with the Tychonov product topology. 


Now let $A$ be a closed discrete subspace of $\IR$. The operators $\partial^D_a$, $a\in A$, compose a bicontinuous operator
$$\partial^D_A:C^\infty(\IR)\to\boldprod^A\IR^D,\;\;\partial^D_A:f\mapsto
(\partial^D_a f)_{a\in A}=(f^{(n)}(a))^{n\in\D}_{a\in A}.$$
Here $\boldprod^A\IR^D$ is a bitopological space whose weak topology is the Tychonov product topology while the strong topology is the box-topology of the space $\square^A\IR^D$.

Now we shall try to locate the image $\partial^D_A(\D^\infty_+(\IR))$ of $\D^\infty_+(\IR)$ in $\boldprod^A\IR^D$. For the singleton $D=\bs 1=\{0\}\subset\w$, this is easy: the image $\partial^{\bs 1}_A(\D^\infty_+(\IR))=\{(f(a))_{a\in A}\in\boldprod^A\IR:f\in\D^\infty_+(\IR)\}$ coincides with the set $\ID^{\bs 1}_+(A)$ of all sequences $(f^0_a)_{a\in A}$ such that
\begin{itemize}
\item $f^0_a<f^0_b$ for any $a<b$ in $A$;
\item  $\inf_{a\in A}f^0_a=-\infty$ if $\inf A=-\infty$, and 
\item $\sup_{a\in A}f^0_a=+\infty$ if $\sup A=+\infty$.
\end{itemize}

For a subset $D\subset\w$ containing $0$ and $1$ we put
$$\ID^D_+(A)=\{(f^n_a)^{n\in D}_{a\in A}\in\boldprod^A\IR^D:(f^0_a)_{a\in A}\in\ID^{\bs 1}_+(A)\mbox{ and }(f^1_a)_{a\in A}\in \IR_+^A\}$$and observe that
$\partial^D_A(\D^\infty_+(\IR))\subset\ID^D_+(A)$. 

The proof of Theorem~\ref{main} will heavily rely on the fact that the differential operator $\partial^\w_\IZ:\D^\infty_+(\IR)\to\ID^\w_+(\IZ)$ has a bicontinuous section $\E_\IZ:\ID^\w_+(\IZ)\to\D^\infty_+(\IR)$.
Given a double sequence $f=(f_a^n)^{n\in\w}_{a\in \IZ}\in\ID^\w_+(\IZ)\subset \boldprod^\IZ\IR^D$ the operator $\E_\IZ$ returns a smooth function $E_\IZ f\in\D_+^\infty(\IR)$ with derivatives $(\E_\IZ f)^{(n)}(a)=f^n_a$, $n\in\w$, at each point $a\in \IZ$.

The construction of the section $\E_\IZ$ will rely on the sections  $\E_{\bold 1}$ and $\E_{\bold 2}$ of the operators $\partial^\w_{\bold 1}$ and $\partial^\w_{\bold 2}$ for the singleton $\bold 1=\{0\}\subset\IR$ and the doubleton $\bold 2=\{0,1\}\subset\IR$.
 
By a {\em section} of a function $f:X\to Y$ we understand any function $g:Y\to X$ such that $f\circ g=\id_Y$.

\section{A section $\E_{\bold 1}$ of the operator $\partial^\w_{\bold 1}$}

In this section we shall construct a weakly continuous section $\E_{\bold 1}$ of the differential operator $$\partial^\w_{\bold 1}:\D^\infty_+(\IR)\to \ID^\w_+(\bold 1)\subset\IR^\w,\;\;\partial^\w_{\bold 1}:f\mapsto (f^{(n)}(0))^{n\in\w}$$ for the singleton $\bold 1=\{0\}\subset \IR$.  Here $\ID^\w_+(\bs 1)=\{(f^n_0)^{n\in\w}\in\IR^\w:f^1_0>0\}\subset\IR^\w$ is a bitopological space whose weak and strong topologies coincide with the Tychonov product topology inherited from $\IR^\w$.

\begin{proposition}\label{sec1}
The operator $\partial^\w_{\bold 1}:\D^\infty_+(\IR)\to\ID^\w_+(\bold 1)$ has a weakly continuous section $\E_{\bold 1}:\ID^\w_+(\bold 1)\to\D^\infty_+(\IR)$. Consequently, $\partial^\w_{\bs 1}(\D^\infty_+(\IR))=\ID^\w_+(\bs 1)$.
\end{proposition}

The proof of this proposition is rather long and requires a lot of preliminary work.

Choose a smooth function $\gamma:\IR\to[0,1]$ such that $\gamma(x)=1$ if $|x|\le1/2$ and $\gamma(x)=0$ if $|x|\ge 1$. 

\begin{lemma}\label{lem_cut} {\rm (1)} The map ${\Bbb R} \ni c \longmapsto \gamma(cx) \in C^\infty({\Bbb R})$  is weakly continuous.  
\vskip 1mm 
{\rm (2)} For any $c > 0$ and $n, s \in \w$ \ we have \hspace{4mm} 
$\ds \left|\big(x^n \gamma(c x) \big)^{(s)}\right| \ \leq \ (2n)^s c^{s-n} \|\gamma\|_s$ \ \ $(x \in {\Bbb R})$ 
\end{lemma}

\begin{proof} (2) $\ds \big(x^n \gamma(c x) \big)^{(s)} = \ 
\sum_{k=0}^s \binom{s}{k}  \, (x^n)^{(k)}\, \gamma(cx)^{(s-k)} 
\ = \sum_{k=0}^{\min \{n,s \}} \hspace{-1mm} \binom{s}{k}\frac{n!}{(n-k)!} x^{n-k}\,\gamma^{(s-k)}(cx)\, c^{s-k}$. \\[1mm]
Here, as expected, $\ds \binom{s}{k}=\frac{s!}{k!(s-k)!}$ stands for the binomial coefficient. 
\vskip 2mm 
Therefore, we have 
\vspace{-2mm} 
$$\begin{array}[t]{lll}
\ds \left|\big(x^n \gamma(c x) \big)^{(s)}\right| & \leq & \ds 
\sum_{k=0}^{\min \{n,s \}} \binom{s}{k}\frac{n!}{(n-k)!}\, |x|^{n-k}\,|\gamma^{(s-k)}(cx)|\, c^{s-k} \\[4mm] 
& \leq & \ds 
\sum_{k=0}^{\min \{n,s \}}\binom{s}{k} \  \frac{n!}{(n-k)!}\, c^{-(n-k)}\,\|\gamma^{(s-k)}\|_0 \, c^{s-k} \\[4mm]  
& \leq & \ds 
c^{s-n} \sum_{k=0}^{\min \{n,s \}} \binom{s}{k} \ n^k \, \|\gamma^{(s-k)}\|_0 
\ \leq \ 
c^{s-n} \, n^s\, \|\gamma\|_s \sum_{k=0}^s \binom{s}{k} \ = \ 
(2n)^s\, c^{s-n}\, \|\gamma\|_s.
\end{array}$$ 
\vskip -4mm 
\end{proof} 

Let $n\in\w$. 
For $b \in {\Bbb R}$ and $c > 0$, consider the function 
$\ds g_n(x) = \frac{b}{n!}\, x^n \gamma(c x) \ \in C^\infty({\Bbb R})$. 

\begin{lemma}\label{lem_n_func} 
{\rm (1)} The map \ ${\Bbb R} \times \IR_+ \ni (b, c) \longmapsto g_n \in C^\infty({\Bbb R})$ \ is weakly continuous. 

{\rm (2)} Suppose $s \geq 0$ and $n \geq s+1$.  
\begin{itemize}
\item[(i)\,] If $c \geq 1$, then \ $\ds \left\| g_n^{(s)} \right\|_0 \ \leq \ |b| (2n)^s c^{-1} \|\gamma\|_n$. 
\vskip 2mm 
\item[(ii)] If $d > 0$ \ and \ $\ds c = n \, 2^{2n} \frac{| b |}{d} \| \gamma \|_n + 3$, \ then \ \ 
$\ds \left\| g_n^{(s)} \right\|_0 \ \leq \ \frac{n^{s-1} d}{\,2^n\,}.$ 
\end{itemize} 
\end{lemma} 

\begin{proof} (2) 
Since $\ds g_n^{(s)}(x) = \frac{b}{n!}\,  \big(x^n \gamma(c x) \big)^{(s)}$\!\!, \,from Lemma~\ref{lem_cut} it follows that 
\[ \left| g_n^{(s)}(x) \right| 
= \frac{|b|}{n!} \left|\big(x^n \gamma(c x) \big)^{(s)}\right| 
\leq |b| \ (2n)^s c^{s-n} \|\gamma\|_s \leq |b| \ (2n)^s c^{s-n} \|\gamma\|_n. \] 
\begin{itemize}
\item[(i)\,] Since $c \geq 1$ and $n \geq s+1$, we have \hspace{5mm} 
$\ds \left| g_n^{(s)}(x) \right| 
\leq |b| (2n)^s c^{-1} \|\gamma\|_n \hspace{5mm} (x \in {\Bbb R}). $ 

\item[(ii)] If $b = 0$, then the assertion is trivial. If $b \neq 0$, 
since $c \geq 1$ and $\ds c > n 2^{2n} \frac{| b |}{d} \| \gamma \|_n$, by (i) 
we have 
\[\ds \left| g_n^{(s)}(x) \right| 
\ \leq \ \frac{|b| (2n)^s \|\gamma\|_n}{c}   
\ \leq \  \frac{|b| (2n)^{s}\| \gamma \|_n}{n 2^{2n}  \frac{| b |}{d} \| \gamma \|_n}
\ = \ \frac{2^s}{2^n}\frac{n^{s-1} d}{\,2^n\,}
\ \leq \ \frac{n^{s-1} d}{\,2^n\,}
\hspace{5mm} (x \in {\Bbb R}).  \] 
\end{itemize}
\vskip -5mm 
\end{proof} 

Let us recall that $\ID^\w_+(\bs 1)=\{(f^n_0)^{n\in\w}\in\IR^\w:f^1_0>0\}$.
Given a sequence $f=(f^n_0)^{n\in\w} \in \ID^\w_+(\bs 1)$, consider the sequence of functions $g_n\in C^\infty(\IR)$ defined by  
$$g_0(x) = f^0_0,\;\;g_1(x) = f^1_0\cdot x, \mbox{ and } 
g_n(x) = \frac{f^n_0}{n!}\, x^n \gamma(c_n x), \mbox{where }c_n = n \, 2^{2n} \frac{| f^n_0 |}{f^1_0} \| \gamma \|_n + 3\mbox{ for }n\ge 2.$$

The following lemma implies Proposition~\ref{sec1}.  

\begin{lemma} \label{lem_series} 
{\rm (1)} The series \ $\ds \sum_{n=0}^\infty g_n^{(s)}$ converges uniformly 
for any $s\in\w$. 

{\rm (2)} The function $\ds g := \sum_{n=0}^\infty g_n$ satisfies the following  conditions: 
\vspace{-1mm} 
\begin{itemize} 
\item[]  \ {\rm (i)} \ $g \in C^\infty({\Bbb R})$, \ \ 
{\rm (ii)} \ $\ds g^{(s)} 
= \sum_{n=0}^\infty g_n^{(s)}$ \ $(s\in\w)$, \ \ 
{\rm (iii)} \ $\partial_{\bold 1}^\w g= f$, \ \ 
\vskip 2mm 
\item[] 
{\rm (iv)} \ $g'(x) > 0$, \ \ 
{\rm (v)} \ $g(x) = f^0_0 + f^1_0\cdot x$ \ \ $(|x| \geq 1/3)$ \ \ and \ \ 
{\rm (vi)} \ $g\in {\mathcal D}^\infty_+({\Bbb R})$. 
\end{itemize}
\vskip 1mm

{\rm (3)} The map \ $\E_{\bold 1} : \ID^\w_+(\bs 1) \to {\mathcal D}^\infty_+({\Bbb R}) \subset C^\infty({\Bbb R})$, \ $\E_{\bold 1}:f\mapsto g$, \ is weakly continuous.
\vskip 2mm

{\rm (4)} $\partial^\w_{\bs 1}(\D^\infty_+(\IR))=\ID^\w_+(\bs 1)$. 
\end{lemma} 

\begin{proof} By Lemma~\ref{lem_n_func}\,(2-ii), 
for $s\ge 0$ and $n\ge s+2$ we have \  
$\ds \left\| g_n^{(s)} \right\|_0 \ \leq \ \frac{n^{s-1} f^1_0}{\,2^n\,}$. 
\smallskip

(1) Since $\ds \sum_{n=0}^\infty \frac{n^s}{2^n} < \infty$,  we have 
$\ds \sum_{n=2}^\infty \left\| g_n^{(s)} \right\|_0 < \infty$, 
which implies the assertion.
\smallskip

(2) The assertions (i), (ii) follow from (1) and the basic properties of uniformly convergent series. 
\vskip 2mm 
\begin{enumerate}
\item[(iii)]  
Since $g_0(x) = f^0_0$, \ $g_1(x) = f^1_0\cdot x$, \ 
$\ds g_n(x) = \frac{f^n_0}{n!}\,x^n$ \ ($|x| < \frac1{2c_n}$) \ ($n \geq 2$), 
it follows that \\[1mm]  
$\ds g_n^{(s)}(0) = f^n_0 \,\delta_{n,s}$ \ $(n \geq 0)$ \ \ and \ \ 
$\ds g^{(s)}(0) = \sum_{n=0}^\infty g_n^{(s)}(0) = f^s_0$. 
\vskip 1mm
\item[(iv)] Since $\ds g'(x) = \sum_{n=0}^\infty g_n^{(1)}(x)
= f^1_0 + \sum_{n=2}^\infty g_n^{(1)}(x)$, it follows that 
$$\ds |g'(x) - f^1_0| 
\ = \ \left| \sum_{n=2}^\infty  g_n^{(1)}(x) \right| 
\ \leq \ \sum_{n=2}^\infty \big|  g_n^{(1)}(x) \big| 
\ \leq \ \sum_{n=2}^\infty \big\|  g_n^{(1)} \big\|_0 
\ \leq \ \sum_{n=2}^\infty \frac{f^1_0}{2^n} 
\ = \ \frac{f^1_0}{2},$$ 
and so \ $g'(x) \ \geq \ \frac{1}{2}f^1_0 \ > \ 0$.  
\vskip 1mm 
\item[(v)\,] The conclusion follows from the equality $g_n(x) = 0$ holding for all $n \geq 2$ and $|x| \geq 1/3$. 
\item[(vi)] By (v) $g$ is surjective. Thus, by (i) and (iv) we have $g \in {\mathcal D}_+^\infty(\IR)$.  
\end{enumerate}
\smallskip

(3) The weak continuity of the map $\E_{\bold 1}:\ID^\w_+(1)\to\D^\infty_+(\IR)$, $\E_{\bold 1}:f\mapsto g$, will follow as soon as we check that for every $s\in\w$ the map  
$$\partial^{(s)}\circ \E_{\bold 1}: \ID^\w_+(1)\to C(\IR),\quad \partial^{(s)}\circ \E_{\bold 1}:f\mapsto g^{(s)},$$ is weakly continuous.
   
By (2)(ii) we have $g^{(s)} = \sum_{n=0}^\infty g_n^{(s)}$. 
Note that the following maps are weakly continuous: 
\[ \mbox{$\ID^\w_+(1) \ni f \longmapsto g_n \in C^\infty({\Bbb R})$, \hspace{5mm} 
$\ID^\w_+(1) \ni f \longmapsto g_n^{(s)} \in C({\Bbb R})$.} \]  
The weak continuity of $g_n$ in $f$ follows from the definition of $g_n$ and  
the weak continuity of $g_n^{(s)}$ in $f$ follows from the weak continuity of the map 
$C^\infty({\Bbb R}) \ni h \longmapsto h^{(s)} \in C({\Bbb R})$. 
Since 
$\ds \left\| g_n^{(s)} \right\|_0 \ \leq \ \frac{n^{s-1} f^1_0}{\,2^n\,}$ \  
($n \geq s + 2$), we see that the map $\partial^{(s)}\circ \E_{\bold 1}$ is weakly continuous. 
\smallskip

(4) The equality $\partial^\w_{\bs 1}(\D^\infty_+(\IR))=\ID^\w_+(1)$ follows from (2-iii).
\end{proof}

\section{A section $\E_{\bold 2}$ of the operator $\partial^\w_{\bold 2}$}

In this section for the doubleton $\bs 2=\{0,1\}\subset\IR$ we construct a continuous section $\E_{\bs 2}$ of the differential operator 
$$\partial^\w_{\bold 2}:\D^\infty_+(\IR)\to \ID^\w_+(\bold 2),\quad\partial^\w_{\bold 2}:f\mapsto(f^{(n)}(a))^{n\in\w}_{a\in\bold 2}.$$ Here $$\ID^\w_+(\bs 2)=\{(f^n_a)^{n\in\w}_{a\in \bs 2}:f^0_0<f^0_1,\;\;(f^1_a)_{a\in\bs 2}\in\IR_+^{\bs 2}\}\subset\boldprod_{a\in\bs 2}\IR^\w$$ is a bitopological space whose weak and strong topologies coincide with the Tychonov product topology inherited from $\boldprod_{a\in\bs 2}\IR^\w$. 

\begin{proposition}\label{sec2}
The operator $\partial^\w_{\bold 2}:\D^\infty_+(\IR)\to\ID^\w_+(\bold 2)$ has a  weakly continuous section $\E_{\bold 2}:\ID^\w_+(\bold 2)\to\D^\infty_+(\IR)$. Consequently, $\partial^\w_{\bs 2}(\D^\infty_+(\IR))=\ID^\w_+(\bs 2)$.
\end{proposition} 

The proof of this proposition requires some preparatory work.
First, let us prove a lemma that will allow to join smooth increasing functions.

\begin{lemma}\label{lem_join} If $f, g \in C^r_+(\IR)$, $f \le g$, and 
$\alpha \in C^\infty(\IR, [0,1])$, $\alpha' \geq 0$, then 
$h := (1 - \alpha) f + \alpha g \in C^r_+(\IR).$ 
\end{lemma} 

\begin{proof} For $r \geq 1$ we have $h' = (-\alpha') f + (1-\alpha) f' + \alpha' g + \alpha g'  
= \alpha' (g - f) + (1 - \alpha) f' + \alpha g' > 0.$

If $r=0$, then for any real numbers $x<y$ we get $\alpha(x)\le\alpha(y)$ and 
$$
\begin{aligned}
h(x)&=(1-\alpha(x))\,f(x)+\alpha(x)\,g(x)<(1-\alpha(x))\,f(y)+\alpha(x)\,g(y)=\\
&=f(y)+\alpha(x)\,(g(y)-f(y))\le f(y)+\alpha(y)\,(g(y)-f(y))=h(y),
\end{aligned}$$
witnessing that $h\in C^0_+(\IR)$.
\end{proof} 

\begin{lemma}\label{leftE} There exists a weakly continuous map $\E_0:\ID^\w_+(\bold 2)\to C^\infty_+(\IR)$ such that for every $f=(f^n_a)^{n\in\w}_{a\in\bs 2}\in\ID^\w_+(\bold 2)$ the function $\E_0 f$ has the following properties:
\begin{enumerate}
\item $(\E_0 f)^{(n)}(0)=f^n_0$ for all $n\in\w$;
\item $\E_0f(\IR)=\big(-\infty,(f^0_0+f^0_1)/2\big)$.
\end{enumerate}
\end{lemma}

\begin{proof} Choose any smooth function $\alpha:\IR\to[0,1]$ such that $\alpha'\ge0$, $(-\infty,\frac13]\subset \alpha^{-1}(0)$ and $[\frac23,+\infty)\subset\alpha^{-1}(1)$.
Also fix a smooth function $\beta\in C^\infty_+(\IR)$ such that $\beta(\IR)=(\tfrac13,\tfrac12)$.

By Proposition~\ref{sec1}, the differential operator $$\partial^\w_{\bold 1}:\D^\infty_+(\IR)\to\ID^\w_+(\bold 1),\;\partial^\w_{\bold 1}:f\mapsto (f^{(n)}(0))^{n\in\w},$$ has a weakly continuous section $\E_{\bold 1}:\ID^\w_+(\bold 1)\to\D^\infty_+(\IR)$.

Given a double sequence $f=(f^n_a)^{n\in\w}_{a\in\bs 2}\in\ID^\w_+(\bold 2)$, consider the sequence $f_0=(f_0^n)^{n\in\w}\in\ID^\w_+(\bold 1)$ and the smooth function $\E_{\bold 1}f_0\in\D^\infty_+(\IR)$. Since $\E_{\bs 1}f_0(0)=f_0^0<f^0_1$, the number $c_f=(\E_{\bs 1}f_0)^{-1}(\frac23f^0_0+\frac13f^0_1)$ is positive. It follows that for every $x\in[0,c_f]$ we get 
$\E_{\bs 1}f_0(x)\le \frac23f^0_0+\frac13f^0_1<\beta_f(x)$, where 
$$\beta_f(x)=f^0_0+(f^0_1-f^0_0)\cdot\beta(x)\in \mbox{$\big(\frac23f^0_0+\frac13f^0_1,\frac12f^0_0+\frac12f^0_1\big).$}$$

 By Lemma~\ref{lem_join}, the function $$\E_0 f:x\mapsto (1-\alpha(c_fx))\cdot\E_{\bs 1}f_0(x)+\alpha(c_fx)\cdot \beta_f(x)$$
belongs to $C^\infty_+(\IR)$. The choice of the functions $\alpha$ guarantees that $\E_0 f(x)=\beta_f(x)$ if $x\ge \frac23c_f$ and $\E_0 f(x)=\E_{\bs 1}f_0(x)$ if  $x\le\frac13 c_f$. Consequently,
$$\E_0f(\IR)=(-\infty,\sup_{x\in\IR}\beta_f(x))=\big(-\infty,\tfrac{f^0_0+f^0_1}2\big)$$and
 $$(\E_0 f)^{(n)}(0)=(\E_{\bold 1}f_0)^{(n)}(0)=f^n_0$$ for all $n\in\w$.

The construction of the function $\E_0 f$ implies that the map
$$\E_0:\ID^\w_+(\bold 2)\to C^\infty_+(\IR),\;\;\E_0:f\mapsto\E_0 f$$is weakly continuous.
\end{proof}

By analogy we can prove

\begin{lemma}\label{rightE} There exists a weakly continuous map $\E_1:\ID^\w_+(\bold 2)\to C^\infty_+(\IR)$ such that for every $f=(f^n_a)^{n\in\w}_{a\in\bs 2}\in\ID^\w_+(\bold 2)$ the function $\E_1 f$ has the following properties:
\begin{enumerate}
\item $(\E_1 f)^{(n)}(1)=f^n_1$ for all $n\in\w$;
\item $\E_1f\,(\IR)=\big(\frac{f^0_0+f^0_1}2,+\infty)$.
\end{enumerate}
\end{lemma}

\begin{proof}[\bf Proof of Proposition~\ref{sec2}]  
Let \ $\E_0,\E_1:\ID^\w_+(\bold 2)\to C^\infty_+({\Bbb R})$ be the weakly continuous maps given by Lemmas~\ref{leftE} and \ref{rightE}. 

Fix a smooth function $\alpha\in C^\infty(\IR,[0,1])$ such that $\alpha'\ge0$, $(-\infty,\frac13]\subset\alpha^{-1}(0)$, and $[\frac23,+\infty)\subset\alpha^{-1}(1)$. 

Given a double sequence $f = (f^n_a)^{n\in\w}_{a\in \bs 2}\in\ID^\w_+(\bs 2)$, consider the map 
$$\E_{\bold 2}f:x\mapsto (1-\alpha(x))\cdot\E_0 f(x)+\alpha(x)\cdot\E_1 f(x),$$ and observe that
$$\E_{\bold 2}f(x)=\begin{cases}
\E_0 f(x)&\mbox{if $x\le \frac13$},\\
\E_1 f(x)&\mbox{if $x\ge \frac23$}
\end{cases}
$$Consequently, $(\E_{\bold 2} f)^{(n)}(0)=(\E_0 f)^{(n)}(0)=f^n_0$ and $(\E_{\bold 2} f)^{(n)}(1)=(\E_1 f)^{(n)}(1)=f^n_1$, witnessing that $\partial^\w_{\bold 2}(\E_{\bold 2}f)=f$.

It follows from $\sup_{x\in\IR}\E_0 f\,(x)=\frac{f^0_0+f^0_1}2=\inf_{x\in\IR}\E_1f\,(x)$ and Lemma~\ref{lem_join} that the map $$\E_{\bold 2}f=\alpha\,\E_0f+(1-\alpha)\E_1f$$ belongs to $C_+^\infty(\IR)$. Being a surjective smooth function with positive derivative, the function $\E_{\bold 2}f$ belongs to the diffeomorphism group $\D^\infty_+(\IR)$.

The weak continuity of the operators $\E_0$ and $\E_1$ imply the weak continuity of the extension operator $$\E_{\bold 2}:\ID^\w_+(\bold 2)\to\D^\infty_+(\IR),\;\E_{\bold 2}:f\mapsto \E_{\bold 2}f.$$
\end{proof}

\section{A section $\E_{\IZ}$ of the operator $\partial^\w_{\IZ}$}

In this section we shall construct a bicontinuous section $\E_\IZ$ of the differential operator 
$$\partial^\w_\IZ:\D^\infty_+(\IR)\to\ID^\w_+(\IZ),\quad \partial^\w_\IZ:f\mapsto(f^{(n)}(a))^{n\in\w}_{a\in\IZ}.$$
Here $$\ID^\w_+(\IZ)=\{(f^n_a)^{n\in\w}_{a\in\IZ}\in\boldprod_{a\in\IZ}\IR^\w:(f^0_a)_{a\in\IZ}\in\ID^{\bold 1}_+(\IZ),\;(f^1_a)_{a\in\IZ}\in\IR_+^\IZ\},$$where$$\ID^{\bs 1}_+(\IZ)=\{(f^0_a)_{a\in \IZ}\in\boldprod_{a\in \IZ}\IR: f^0_a<f^0_{a{+}1}\;(\forall a\in\IZ),\;\;\inf_{a\in\IZ}f^0_a=-\infty,\;\;\sup_{a\in\IZ}f^0_a=+\infty\}.$$
The space $\ID^{\w}_+(\IZ)$ is a subspace of the bitopological space 
$\boldprod^{\IZ}\IR^\w$ whose weak topology is the Tychonov product topology while the strong topology is the topology of the box-product $\square^\IZ\IR^\w$.
The same concerns the subspace $\ID^{\bs 1}_+(\IZ)$ of the bitopological space 
$\boldprod^{\IZ}\IR$ endowed with the Tychonov and box-product topologies.

It is clear that the image $\partial^{\bs 1}_\IZ(\D^\infty_c(\IR))$ of the group $\D^\infty_c(\IR)$ of compactly supported $C^\infty$ diffeomorphisms of the real line lies in the subset
$$\ID^{\bs 1}_c(\IZ)=\{(f^0_a)_{a\in\IZ}\in\ID^{\bs 1}_+(\IZ):|\{a\in\IZ: f^0_a\ne a\}|<\aleph_0\}\subset\ID^{\bs 1}_+(\IZ).$$ 

By analogy, the image $\partial^\w_\IZ(\D^\infty_c(\IR))$ of $\D^\infty_c(\IR)$ lies in the subset
$$\ID^\w_c(\IZ)=\{(f^n_a)^{n\in\w}_{a\in\IZ}\in\ID^\w_+(\IZ):|\{a\in\IZ: (f^n_a)^{n\in\w}\ne \partial^\w_a\id_\IR\}|<\aleph_0\}\subset\ID^\w_+(\IZ).$$ 


\begin{proposition}\label{secZ} The differential operator $\partial^\w_\IZ:\D^\infty_+(\IR)\to\ID^\w_+(\IZ)$ has a bicontinuous section $\E_\IZ:\ID^\w_+(\IZ)\to\D^\infty_+(\IR)$ such that  $\E_\IZ(\ID^\w_c(\IZ))\subset\D^\infty_c(\IR)$. Consequently, $\partial^\w_\IZ(\D^\infty_+(\IR))=\ID^\w_+(\IZ)$ and $\partial^\w_\IZ(\D^\infty_c(\IR))=\ID^\w_c(\IZ)$.
\end{proposition}
 
We shall prove this proposition after some preparatory work.

Given a point $a\in\IR$ and two function $f,g:\IR\to\IR$ with $f(a)=g(a)$ let $f\cup_a g$ be the function equal to $f$ on the ray $(-\infty,a]$ and to $g$ on the ray $[a,+\infty)$. It is clear that $f\cup_a g$ is continuous if so are the functions $f$ and $g$.

In order to describe the smoothness properties of  $f\cup_a g$, for every $n \le r \le \infty$ consider the differential operator 
$$\partial^{\leqslant n}_a:C^r(\IR)\to\prod_{k=0}^n\IR,\;\;\partial^{\leqslant n}_a:f\mapsto(f^{(k)}(a))_{k=0}^n.$$

The next lemma follows from the definition of the derivative and 
the induction on $n\in\w$. 

\begin{lemma}\label{lem_paste} 
If $f,g \in C^n({\Bbb R})$ and $\partial^{\leqslant n}_a f=\partial^{\leqslant n}_a g$, then \ $h  = f \cup_{a} g \ \in C^n({\Bbb R})$ \ and \ 
$h^{(n)} = f^{(n)} \cup_{a} g^{(n)}$. 
\end{lemma} 

\begin{lemma}\label{lem_derivative} 
Suppose $f, g, \phi, \psi \in {\mathcal D}_+^n({\Bbb R})$, and 
let $a \in {\Bbb R}$ and $b=f(a)$.  
\begin{itemize}
\item[(1)] If $\partial_a^{\leqslant n} f = \partial_a^{\leqslant n} g$ and 
$\partial_b^{\leqslant n}\phi = \partial_b^{\leqslant n} \psi$, then 
$\partial_a^{\leqslant n} (\phi \circ f) = \partial_a^{\leqslant n} (\psi \circ g)$. 
\item[(2)] If $\partial_a^{\leqslant n} f = \partial_a^{\leqslant n} g$, then $\partial_b^{\leqslant n} f^{-1} = \partial_b^{\leqslant n} g^{-1}$. In particular, 
\[ \mbox{{\rm (i)} $\partial_a^{\leqslant n} f = \partial_a^{\leqslant n} \id_\IR$ \ iff \ 
$\partial_a^{\leqslant n} f^{-1} = \partial_a^{\leqslant n} \id_\IR$, \hspace{2mm} and \hspace{2mm} 
{\rm (ii)} $\partial_a^{\leqslant n} (\varphi \circ f) = \partial_a^{\leqslant n} \varphi$ \ iff \ 
$\partial_a^{\leqslant n} f = \partial_a^{\leqslant n} \id_\IR $.} \] 
\end{itemize}
\end{lemma}

\begin{proof} 
(1) Lemma~\ref{lem_paste} implies
$$f\cup_{a} g,\;\phi \cup_{b} \psi \in {\mathcal D}_+^n({\Bbb R})\mbox{  and }   
h := (\phi \circ f) \cup_{a} (\psi \circ g) = (\psi \cup_{b} \psi) \circ (f \cup_{a} g) \in {\mathcal D}_+^n({\Bbb R})$$ and hence 
$\partial_a^{\leqslant n} (\phi \circ f) = \partial_a^{\leqslant n} h = \partial_a^{\leqslant n} (\psi \circ g).$ 

(2) Since $f\cup_{a} g \in {\mathcal D}_+^n({\Bbb R})$, we get  
$h := f^{-1} \cup_{b} g^{-1} = (f \cup_{a} g)^{-1} \in {\mathcal D}_+^n({\Bbb R})$ and hence $\partial_b^{\leqslant n} f^{-1} = \partial_b^{\leqslant n} h = \partial_b^{\leqslant n} g^{-1}$. 
\end{proof} 

We shall use Lemma~\ref{lem_derivative} to prove the following self-generalization of Proposition~\ref{sec2}. 

\begin{proposition}\label{sec2a} For any doubleton $A=\{a,b\}\subset\IR$ the differential operator
$$\partial^\w_{A}:\D^\infty_+(\IR)\to \ID^\w_+(A),\quad\partial^\w_A:g\mapsto \big(g^{(n)}(x)\big)^{n\in\w}_{x\in A},$$
admits a weakly continuous section $\E_A:\ID^\w_+(A)\to\D^\infty_+(\IR)$ such that $\E_A\circ \partial^\w_A(\id_\IR)=\id_\IR$.
\end{proposition}

\begin{proof} According to Proposition~\ref{sec2}, the differential operator $\partial^\w_{\bold 2}:\D^\infty_+(\IR)\to \ID^\w_+(\bold 2)$ has a weakly continuous section $\E_{\bold 2}:\ID^\w_+(\bold 2)\to\D^\infty_+(\IR)$.

Write $A=\{a,a+c\}$ for some $c>0$ and consider the affine map $\ell(x)=a+cx$ of the real line which maps the doubleton $\bold 2=\{0,1\}$ onto the doubleton $A$. This map induces two weak homeomorphisms:
$$\mathsf L:\D^\infty_+(\IR)\to \D^\infty_+(\IR),\quad\mathsf L:f\mapsto f\circ \ell,$$
and $$\Lambda:\ID^\w_+(A)\to\ID^\w_+(\bold 2),\;\Lambda:(f^n_a)^{n\in\w}_{a\in A} \mapsto \big(c^n\cdot f^n_{\ell(a)})^{n\in\w}_{a\in\bs 2}.$$ 

For every $f\in C^\infty(\IR)$ and $n\in\w$ we get $(f \circ \ell)^{(n)}(x) = c^n f^{(n)}(\ell(x))$. 
%
%
This yields the commutativity of the following diagram:
$$\begin{CD}
\D^\infty_+(\IR)@>{\mathsf L}>>\D^\infty_+(\IR)\\
@V{\partial^\w_A}VV@VV{\partial^\w_{\bold 2}}V\\
\ID^\w_+(A)@>{\Lambda}>>\ID^\w_+(\bold 2)
\end{CD}
$$

It follows that the map $\tilde\E_A=\mathsf L^{-1}\circ \E_{\bold 2}\circ\Lambda$ is a weakly continuous section of the operator $\partial^\w_A$. However this section can map the double sequence $\partial^\w_A\id_\IR \in \ID^\w_+(A)$ onto some function $g=\tilde\E_A(\partial^\w_A\id_\IR)$ that is not equal to $\id_\IR$. To fix this problem, we modify the map $\tilde\E_A$ to the map $$\E_A:\ID^\w_+(A)\to\D^\infty_+(\IR),\;\;\E_A:f\mapsto (\tilde \E_A f)\circ g^{-1}.$$ This map is weakly continuous because $\D^\infty_+(\IR)$ endowed with the compact-open $C^\infty$ is a topological group. It is clear that $$\E_A\partial^\w_A\id_\IR=(\tilde \E_A\partial^\w_A\id_\IR)\circ g^{-1}=g\circ g^{-1}=\id_\IR.$$

According to Lemma~\ref{lem_derivative}, the equality $\partial^\w_A g=\partial^\w_A\id_A$ implies that $\partial^\w_A g^{-1}=\partial^\w_A\id_\IR$. 
Then for every $f\in\ID^\w_+(A)$ we get 
$$\partial^\w_A\E_A f= \partial^\w_A \Big(\big(\tilde \E_Af\big)\circ g^{-1} \Big)=\partial^\w_A\tilde \E_A f=f,$$ 
witnessing that $\E_A:\ID^\w_+(A)\to\D^\infty_+(\IR)$ is a section of the operator $\partial^\w_A$.
\end{proof}

\begin{proof}[{\bf Proof of Proposition~\ref{secZ}}] We need to construct a bicontinuous section $\E_\IZ:\ID^\w_+(\IZ)\to\D^\infty_+(\IR)$ of the differential operator $\partial^\w_\IZ:\D^\infty_+(\IR)\to\ID^\w_+(\IZ)$ such that $\E_\IZ(\ID^\w_c(\IZ)) \subset \D^\infty_c(\IR)$.

For every $m\in\IZ$ consider the closed interval $\II_m=[m,m{+}1]$ and its boundary $\partial\II_m=\{m,m+1\}$ in $\IR$. It is clear that the projection
$$\pi_m:\ID^\w_+(\IZ)\to\ID^\w_+(\partial\II_m),\quad
\pi_m:(f^n_a)^{n\in\w}_{a\in\IZ}\mapsto (f^n_a)^{n\in\w}_{a\in\partial\II_m}$$ is bicontinuous.

By Proposition~\ref{sec2a}, for every $m\in\IZ$ the differential operator
$\partial^\w_{\partial\II_m}:\D^\infty(\IR)\to\ID^\w_+(\partial\II_m)$ has a weakly continuous section $\E_{\partial\II_m}:\ID^\w_+(\partial\II_n)\to\D^{\infty}_+(\IR)$ such that $(\E_{\partial\II_m}\circ\partial^\w_{\partial\II_m})\id_\IR=\id_\IR$.

Define the map $\E_\IZ:\ID^\w_+(\IZ)\to\D^\infty_+(\IR)$ assigning to each double sequence $f=(f^n_a)^{n\in\w}_{a\in\IZ}\in\ID^\w_+(\IZ)$ the function $\E_\IZ f\in\D^\infty_+(\IR)$ defined by$$(\E_\IZ f)|\II_m=(\E_{\partial\II_m}\pi_m f)|\II_m \ \mbox{ for every $m\in\IZ$}.$$

Let us show that this definition is correct. For every $m\in\IZ$ consider the smooth function $g_m=\E_{\partial\II_m}(\pi_m f)\in\D^\infty_+(\IR)$. It follows that $g_m^{(n)}(m)=f^n_m$ and $g_m^{(n)}(m{+}1)=f^n_{m{+}1}$ for every $n\in\w$, and hence  $g^{(n)}_m(m)=f^n_m=g^{(n)}_{m-1}(m)$. Consequently, $\partial^\w_{m}g_m=\partial^\w_{m}g_{m-1}$ and then the function $g_{m-1}\cup_mg_m$ is smooth with $\partial^\w_m (g_{m-1}\cup_mg_m)=(f^n_m)^{n\in\w}$.
Since $(\E_\IZ f)|[m{-}1,m{+}1]=(g_{m-1}\cup_mg_m)|[m{-}1,m{+}1]$, we conclude that $\E_\IZ f$ is smooth on the open interval $(m{-}1,m{+}1)$ and $\partial^\w_m\E_\IZ f=(f^n_m)^{n\in\w}$. This implies that $\partial^\w_\IZ\E_\IZ f=f$ and hence $\E_\IZ:\ID^\w_+(\IZ)\to\D^\infty_+(\IR)$ is a required section of the differential operator $\partial^\w_\IZ:\D^\infty_+(\IR)\to\ID^\w_+(\IZ)$.

Using the fact that for every $f\in\ID^\w(\IZ)$ the values of the function $\E_\IZ f$ on each interval $\II_m$, $m\in\IZ$, depend only on the projection $\pi_m f\in\ID^\w_+(\partial\II_m)$, we can show that the map $\E_\IZ:\ID^\w_+(\IZ)\to\D^\infty(\IR)$ is bicontinuous.

Finally we check that $\E_\IZ(\ID^\w_c(\IZ))\subset\D^\infty_c(\IR)$. Take any $f=(f^n_a)^{n\in\w}_{a\in\IZ}\in\ID^\w_c(\IZ)$. By the definition of $\ID^\w_c(\IZ)$, there exists $M\in\w$ such that $(f^n_a)^{n\in\w}=\partial^\w_a\id_\IR$ for every $a\in\IZ$ with $|a|\ge M-1$.

 Then for every $m\in\IZ$ with $|m|> M$ we get $\pi_m f=\partial^\w_{\partial\II_m}\id_\IR$ and $\E_{\partial\II_m}\pi_m f=\id_\IR$. Consequently, $\E_\IZ f|\II_m=\E_{\partial\II_m}\pi_mf|\II_m=\id_\IR|\II_m$ for all $m>|M|$, which means that $\E_\IZ f\in\D^\infty_c(\IR)$.
\end{proof}

\section{The proof of Theorem~\ref{main}}

We need to prove that for every $r \in \w \cup \{ \infty \}$ the triple of diffeomorphism groups $(\D^r(\IR),\D^r_+(\IR),\D^r_c(\IR))$ endowed with the compact-open and Whitney $C^r$ topologies is bihomeomorphic to the triple of homeomorphism groups $(\HH(\IR),\HH_+(\IR),\HH_c(\IR))$ endowed with the compact-open and Whitney topologies.

First we show that the couple $(\D^r_+(\IR),\D^r_c(\IR))$ is bihomeomorphic to $(\HH_+(\IR),\HH_c(\IR))$.

Consider the bicontinuous differential operator $$\partial^{\leqslant r}_\IZ:\D^r_+(\IR)\to\ID^\w_+(\IZ)$$assigning to each $C^r$ diffeomorphism $f\in\D^r_+(\IR)$ the double sequence $(f^n_a)^{n\in\w}_{a\in \IZ}$ where
$$f^n_a=\begin{cases}
f^{(n)}(a)&\mbox{if $n\le r$,}\\
1&\mbox{if $r<n=1$,}\\
0&\mbox{if $r<n\ne 1$}.
\end{cases}
$$
It is clear that $\partial^{\ls r}_\IZ$ maps $\D^r_+(\IZ)$ and $\D^r_c(\IR)$ into the subspaces
$$\ID^{\leqslant r}_+(\IZ)=\{(f^n_a)^{n\in\w}_{a\in\IZ}\in\ID^\w_+(\IZ):f^n_a=\partial^{(n)}_a\id_\IR\mbox{ for all $n>r$ and $a\in\IZ$}\}$$ and
$$\ID^{\ls r}_c(\IZ)=\ID^{\ls r}_+(\IZ)\cap\ID^\w_c(\IZ)$$
of the bitopological space $$\ID^\w_+(\IZ)=\{(f^n_a)^{n\in\w}_{a\in\IZ}:(f^0_a)_{a\in\IZ}\in\ID^{\bs 1}_+(\IZ),\;(f^1_a)_{a\in\IZ}\in\boldprod^\IZ\IR_+\}.$$

By the definition of the spaces $\ID^{\ls r}_+(\IZ)$ and $\ID_c^{\ls r}(\IZ)$, 
the correspondence 
$$(f^n_a)^{n\in\w}_{a\in\IZ} \longleftrightarrow \big((f^0_a)_{a\in\IZ}, (\log f_a^1, (f^n_a)^{n=2, \dots, r})_{a\in\IZ}\big)$$ 
yields the following conclusion: 

\begin{claim}\label{claim1} The pair $(\ID^{\ls r}_+(\IZ),\ID_c^{\ls r}(\IZ))$ is bihomeomorphic to the pair $(\ID_+^{\bs 1}(\IZ)\times\boldprod^\IZ\IR^r,\ID_c^{\bs 1}(\IZ)\times\boldsum^\IZ\IR^r)$.
\end{claim}

Lemma~\ref{lem_derivative} implies that the subset $$\D^r_+(\IR;\IZ)=\{f\in\D^r_+(\IR):\partial^{\ls r}_\IZ f=\partial^{\ls r}_\IZ\id_\IR\}$$ is a (weakly) closed subgroup of $\D^r_+(\IR)$. 

The bitopological structure of the subgroup $\D^r(\IR;\IZ)$ can be described as follows.
Using Lemma~\ref{lem_paste}, one can check that the map $$\Psi:\D^r(\IR;\IZ)\to\boldprod_{m\in\IZ}C^r(\II_m),\;\;\Psi:f\mapsto (f|\II_m)_{m\in\IZ}$$ determines a bihomeomorphism of the group $\D^r(\IR;\IZ)$ onto the product $\boldprod_{m\in\IZ}\D^r_{\partial}(\II_m)$ of the
groups $$\D^r_\partial(\II_m)=\{f\in C^r_+(\II_m):\partial^{\ls r}_{\partial\II_m}f=\partial^{\ls r}_{\partial\II_m} \id_\IR\}\subset C^r(\II_m).$$
The product $\boldprod_{m\in\IZ}\D^r_\partial(\II_m)$ is considered as a bitopological space carrying the Tychonov and box-product topologies.

The bihomeomorphism $\Psi$ maps the subgroup $\D^r_c(\IR;\IZ)=\D^r(\IR;\IZ)\cap\D^r_c(\IR)$ of $\D^r(\IR;\IZ)$ onto the subspace $$\boldsum_{m\in\IZ}\D^r_\partial(\II_m)=\{(f_m)_{m\in\IZ}\in\boldprod_{m\in\IZ}\D^r_\partial(\II_m):|\{m\in\IZ:f_m\ne\id_{\II_m}\}|<\aleph_0\}$$ of 
$\boldprod_{m\in\IZ}\D^r_\partial(\II_m)$.

Observe that for every  integer number $m\in\IZ$  the group $\D^r_\partial(\II_m)$ is a convex subset of the Fr\'echet space $C^r(\II_m)$. Consequently, this group is an absolute retract according to the Borsuk-Dugundji Theorem \cite{Dug}, \cite{Bor}. Being a Polish non-locally compact absolute retract, the group $\D^r_\partial(\II_m)$ is homeomorphic to the separable Hilbert space $l_2$ by the Dobrowolski-Toru\'nczyk Theorem \cite{DT}.

Now we see that the pair $(\boldprod_{m\in\IZ}\D^r_\partial(\II_m),\boldsum_{m\in\IZ}\D^r_\partial(\II_m))$ is bihomeomorphic to the pair $(\boldprod^\IZ l_2,\boldsum^\IZ l_2)$. Taking into account that the former pair is bihomeomorphic to 
$(\D^r(\IR;\IZ),\D^r_c(\IR;\IZ))$, we get

\begin{claim}\label{claim2} The pair $(\D^r(\IR;\IZ),\D^r_c(\IR;\IZ))$ is bihomeomorphic to $(\boldprod^\IZ l_2,\boldsum^\IZ l_2)$.
\end{claim}

By Proposition~\ref{secZ}, the operator $\partial^\w_\IZ:\D^\infty_+(\IR)\to\ID^\w_+(\IZ)$ has a bicontinuous section $\E_\IZ:\ID^\w_+(\IZ)\to\D^\infty_+(\IR)$ such that $\E_\IZ(\ID^\w_c(\IZ))\subset\D^\infty_c(\IR)$. It is easy to see that $\E_\IZ$ restricted to $\ID^{\ls r}_+(\IZ)\subset\ID^\w_+(\IZ)$ is a bicontinuous section of the operator $\partial^{\ls r}_\IZ:\D^r_+(\IR)\to\ID^{\ls r}_+(\IZ)$.

Now consider the map 
$$\Phi:\ID^{\ls r}_+(\IZ)\times\D^r(\IR;\IZ)\to\D^r_+(\IR),\quad \Phi:(f,g)\mapsto (\E_\IZ f)\circ g.$$ 
The bicontinuity of the map $\Phi$ follows from the bicontinuity of the section $\E_\IZ$ and the bicontinuity of the group operation on $\D^r_+(\IR)$.

In fact, the map $\Phi$ is a bihomeomorphism. Its inverse is defined by
$$\Phi' : \D^r_+(\IR) \to \ID^{\ls r}_+(\IZ)\times\D^r(\IR;\IZ),\quad \Phi' :h\mapsto(\partial^{\ls r}_\IZ h,(\E_\IZ\partial^{\ls r}_\IZ h)^{-1}\circ h).$$ 
For every $h\in\D^r_+(\IR)$ and $g=\E_\IZ\partial^{\ls r}_\IZ h$, 
since $\partial^{\ls r}_\IZ g=\partial^{\ls r}_\IZ\E_\IZ\partial^{\ls r}_\IZ f=\partial^{\ls r}_\IZ f$, from Lemma~\ref{lem_derivative} it follows that $\partial^{\ls r}_\IZ g^{-1}=\partial^{\ls r}_\IZ f^{-1}$ and $\partial^{\ls r}_\IZ (g^{-1}\circ f)=\partial^{\ls r}_\IZ (f^{-1}\circ f)=\partial^{\ls r}_\IZ\id_\IR$, 
which witnesses that $g^{-1}\circ f\in\D^r(\IR;\IZ)$. 
The bicontinuity of the map $\Phi'$ follows from the bicontinuity of the maps $\E_\IZ$, $\partial^{\ls r}_\IZ$ and the group operations on $\D^r(\IR)$. 
It also follows from $\E_\IZ(\ID^\w_c(\IZ))\subset\D^\infty_c(\IR)$ that $\Phi'(\D^r_c(\IR))=\ID^{\ls r}_c(\IZ)\times\D^r_c(\IR;\IZ)$.
This means that the pair $(\D^r_+(\IR),\D^r_c(\IR))$ is bihomeomorphic to the pair 
$(\ID^{\ls r}_+(\IZ)\times \D^{r}(\IR;\IZ),\ID^{\ls r}_c(\IZ)\times \D^{r}_c(\IR;\IZ))$. By Claims~\ref{claim1} and \ref{claim2}, the latter pair is bihomeomorphic to the pair 
\[\label{longpair}
(\ID^{\bold 1}_+(\IZ)\times\boldprod^\IZ\IR^r\times\boldprod^\IZ l_2,
\ID^{\bold 1}_c(\IZ)\times\boldsum^\IZ\IR^r\times\boldsum^{\IZ}l_2)\big),\]
which is bihomeomorphic to $$(\ID^{\bold 1}_+(\IZ)\times\boldprod^\IZ l_2,
\ID^{\bold 1}_c(\IZ)\times\boldsum^{\IZ}l_2)\big).$$

Thus we have proved 

\begin{lemma} For every $r \in \w \cup \{ \infty \}$ the pair $(\D^r_+(\IR),\D^r_c(\IR))$ is bihomeomorphic to the pair $\big(\ID^{\bold 1}_+(\IZ)\times\boldprod^\IZ l_2,
\ID^{\bold 1}_c(\IZ)\times\boldsum^{\IZ}l_2)\big)$. 
\end{lemma}

This lemma implies that for every $r \in \w \cup \{ \infty \}$ the pair 
$(\D^r_+(\IR),\D^r_c(\IR))$ is bihomeomorphic to  the pair $(\D^0_+(\IR),\D^0_c(\IR))=(\HH_+(\IR),\HH_c(\IR))$.
Consequently, there is a bihomeomorphism $\mathsf H:\D^r_+(\IR)\to\HH_+(\IR)$ such that $\mathsf H(\D^r_c(\IR))=\HH_c(\IR)$.

It remains to extend the bihomeomorphism $\mathsf H$ to a bihomeomorphism $\tilde{\mathsf H}:\D^r(\IR)\to\HH(\IR)$. For this fix any diffeomorphism $\phi\in\D^\infty(\IR)\setminus D^\infty_+(\IR)$. For example, we can take $\phi(x)=-x$. 

Extend the bihomeomorphism $\mathsf H$ to a bihomeomorphism $\tilde {\mathsf H}:\D^r(\IR)\to\HH(\IR)$ assigning to each $f\in\D^r(\IR)$ the diffeomorphism
$$\tilde {\mathsf H}f=\begin{cases}\mathsf H f&\mbox{if $f\in\D^r_+(\IR)$}\\
\phi^{-1}\circ \mathsf H(f\circ\phi)&\mbox{if $f\in\D^r(\IR)\setminus\D^r_+(\IR)$}.
\end{cases}
$$ It is clear that the bihomeomorphism $\tilde{\mathsf H}$ maps the pair $(\D^r_+(\IR),\D^r_c(\IR))$ onto the pair 
$(\HH_+(\IR),\HH_c(\IR))$. 

\section{Open Problems}

By Corollaries~\ref{cor1} and \ref{cor2}, for every $r \in \w \cup \{ \infty \}$ the pair $(\D^r_+(\IR),\D^r_c(\IR))$ is weakly and strongly homeomorphic to the pair $(\boldprod^\w l_2,\boldsum^\w l_2)$.

\begin{problem} Are the pairs $(\D^r_+(\IR),\D^r_c(\IR))$, $r \in \w \cup \{ \infty \}$,  bihomeomorphic to the pair $(\boldprod^\w l_2,\boldsum^\w l_2)$?
\end{problem}

The answer to this problem is affirmative if the answer to the following problem is affirmative.

\begin{problem} Is the pair $(\ID^{\bs 1}_+(\IZ),\ID^{\bs 1}_c(\IZ))$ bihomeomorphic to the pair $(\boldprod^\w \IR,\boldsum^\w \IR)$?
\end{problem}

By  \cite{BMS}, the pair $(\ID^{\bs 1}_+(\IZ),\ID^{\bs 1}_c(\IZ))$ is weakly and strongly homeomorphic to the pair $(\boldprod^\w \IR,\boldsum^\w \IR)$.
\smallskip

Our next question asks if Theorem~\ref{main} can be generalized to all smooth manifolds.

\begin{problem} Are the triples $(\D^r(M),\D^r_+(M),\D^r_c(M))$ and $(\HH(M),\HH_+(M),\HH_c(M))$ bihomeomorphic for any connected noncompact orientable $C^r$ manifold $M$?
\end{problem}

Here $\D^r(M)$, $\D^r_+(M)$ and $\D^r_c(M)$ are the groups of $C^r$ diffeomorphisms, orientation-preserving $C^r$ diffeomorphisms and 
compactly supported $C^r$ diffeomorphisms of $M$, respectively. 
Those groups are bitopological spaces whose weak topology is the compact-open $C^r$ topology while the strong topology is the Whitney $C^r$ topology (very strong $C^r$ topology in the sense of \cite{Illman}). 
\medskip

Finally, let us discuss the problem of bitopological characterization of spaces of the form $\boldprod^\w X$ or $\boldsum^\w X$ for simple topological spaces $X$ like $\II=[0,1]$, $\II^\w$, $\IR$, or $l_2$. The bitopological spaces $\boldsum^\w\II$ and $\boldsum^\w\II^\w$ have been characterized in \cite{BS}. This characterization implies that the bitopological spaces $\boldsum^\w\II$ and $\boldsum^\w\IR$ are bihomeomorphic.

\begin{problem}\label{prob4} Give a bitopological characterization of the bitopological space $\boldsum^\w l_2$.
\end{problem} 

Using the technique of \cite{Ban} and \cite{BS} it can be shown that the bitoplogical space $\boldsum^\w l_2$ is bihomeomorphic to $l_2\times\boldsum^\w\IR$. So, Problem~\ref{prob4} can be reformulated as the problem of bitopological characterization of the bitopological space $l_2\times \boldsum^\w\IR$.

It should be mentioned that the topological characterizations of the topological spaces  $\mbox{\large$\Sigma$}^\w l_2$ and  $\cbox^\w l^2$ (composing the bitopologocal space $\boldsum^\w l_2$) are known, see \cite{BM}, \cite{BR}.

The bitopological equivalence of the spaces $\boldsum^\w l_2$ and $l_2\times\boldsum^\w\IR$ suggests the following question.

\begin{problem} Are the bitopological spaces $\boldprod^\w l_2$ and $l_2\times\boldprod^\w\IR$ bihomeomorphic?
\end{problem} 

In order to answer this problem it would be helpful to know the answer to the following (apparently, very difficult) problem.

\begin{problem} Give a bitopological characterization of the bitopological spaces  $\boldprod^\w \IR$ and $\boldprod^\w l_2$.
\end{problem}

By Anderson's Theorem \cite{An}, the spaces $\boldprod^\w \IR$ and $\boldprod^\w l_2$ are weakly homeomorphic. On the other hand, these spaces are not strongly homeomorphic  because the connected components of those spaces endowed with the strong topology are homeomorphic to the (non-homeomorphic) spaces $\cbox^\w \IR$ and $\cbox^\w l_2$, respectively. 

The topological characterization of the spaces $\mbox{\large$\Pi$}^\w \IR$ and $\mbox{\large$\Pi$}^\w l_2$ was given by Toru\'nczyk \cite{Tor}. An analogous problem for the respective box-products is open.

\begin{problem} Give a topological characterization of the spaces $\square^\w \IR$ and $\square^\w l_2$.
\end{problem}

\section{Acknowledgment}

The authors thank the referee whose constructive criticism resulted in considerable improvement of the presentation of this paper.

\end{document}